\documentclass{amsart}
\usepackage{amsmath}
\usepackage{amsfonts}
\usepackage{amsthm}
\usepackage{url}
\usepackage{microtype}
\title[Sliced Wasserstein Comparison Constant]{A lower bound on the sliced Wasserstein comparison constant for translated measures}
\author{Avery Noroozi}
\date{}
\address{Department of Mathematics, Purdue University}
\email{noroozi@purdue.edu}
\begin{document}

\begin{abstract}
We provide closed-form expressions for standard 1-Wasserstein distance ($\mathcal{W}_1$) and sliced 1-Wasserstein distance ($\mathcal{SW}_1$) for translated probability measures on a finite dimensional Hilbert space $V$. Furthermore, we construct an explicit lower bound for the constant $c_d$ established by Carlier, Figalli, Mérigot, and Wang between $\mathcal{W}_1$ and $\mathcal{SW}_1$ \cite{CFMW2025} for translated probability measures on $\mathbb{R}^d$. Sliced Wasserstein distances are computationally efficient and used in image generation and other applications, which raises natural questions in computational efficiency versus theoretical correctness. A lower bound on the constant $c_d$ gives information about the cost of that efficiency.
\end{abstract}

\maketitle
\vspace{-0.7cm}
\section{Introduction and Background}
In this text, we will explore two metrics and their ratios for fixed probability measures. Specifically, we consider 1-Wasserstein distance ($\mathcal{W}_1)$ and sliced 1-Wasserstein distance ($\mathcal{SW}_1$), which serves as a computationally cheaper alternative to the standard metric. It projects multi-dimensional distributions onto 1D lines on the unit sphere. 1D optimal transport has a simple form (monotone rearrangement) optimal for all $p$. Sliced 1-Wasserstein distance, achieves $\mathcal{O}(n\text{log}n)$ time complexity, far superior to the time complexity of 1-Wasserstein distance, $\mathcal{O}(n^3\log n)$ at best \cite{LeEtAl2022}.
\\
\\
\indent We offer a concrete example in $\mathbb{R}^3$ shifting the uniform probability measure on the unit disc on the $xy$-plane by $\varepsilon > 0$ in the $z$-direction, computing $\mathcal{W}_1$ and $\mathcal{SW}_1,$ and finally examining their individual values and ratio $\frac{\mathcal{SW}_1}{\mathcal{W}_1},$ giving us a constant value of $\frac{1}{2}$. We then generalize further to finite dimensional Hilbert spaces which we will denote by $V$ unless otherwise specified with $\operatorname{dim}(V) = d$. We find closed-form expressions for $\mathcal{W}_1$ and $\mathcal{SW}_1$ for general translated probability measures in  $\mathcal{P}_1(V)$, the set of probability measures on $V$ with finite first moment.
\\
\\
We consider $\mu, \nu \in \mathcal{P}_1(V)$ such that $T_{\#}\mu = \nu,$ $T: V \rightarrow V$ with $T(x) = x+v$ with $x,\  v \in V, \ v \neq 0.$ We analyze how the ratio $\frac{\mathcal{SW}_1}{\mathcal{W}_1}$ behaves in general dimension $d$ that shows a dependency influenced by the Gamma function $\Gamma.$ Lastly, we provide an explicit lower bound on the constant $c_d$ established by  Carlier, Figalli, Mérigot, and Wang \cite{CFMW2025} for translated measures. For the lower bound, we consider the Euclidean case.
\\
\\
We proceed with some definitions and terminology. The 1-Wasserstein distance between two probability measures $\mu,\nu$ is defined
\[
\mathcal{W}_1(\mu, \nu)= \inf_{\pi \in \Pi(\mu,\nu)} \int_{V \times V} ||x-y|| \, d\pi(x,y),
\]
where $\Pi(\mu,\nu)$ denotes the set of all valid transport plans from $\mu$ to $\nu,$ i.e. joint probability distributions with $\mu,\nu$ as the marginal distributions \cite{Wasserman}. 
\\
\\
Alternatively, the sliced 1-Wasserstein distance is defined,

\[
\mathcal{SW}_1(\mu,\nu)
=
\int_{{\mathbb{S}}^{d-1}}
\mathcal{W}_1\bigl(P_{\theta\#}\mu, P_{\theta\#}\nu\bigr)
\,d\sigma(\theta),
\]
where $d\sigma(\theta)$ is the uniform probability measure on the unit sphere $\mathbb{S}^{d-1} = \{\theta \in V: ||\theta|| = 1\}$ and $P_\theta(v)$ is the inner product $\langle v, \theta \rangle.$ In 2025, Carlier, Figalli, Mérigot, and Wang \cite{CFMW2025} proved the sharp comparison between $\mathcal{W}_1$ and $\mathcal{SW}_1,$ for any dimension $d \geq 2,$ there exists a constant $c_d$ such that for all $R > 0$ and all $\mu,\nu \in \mathcal{P}(\mathbb{R}^d)$ supported on $B_R,$ 
\[
\mathcal{W}_1(\mu,\nu) \leq c_dR^{\frac{d-1}{d}}\mathcal{SW}_1(\mu,\nu)^{\frac{1}{d}} 
\]
where $B_R$ is the ball of radius $R$ centered at 0. 

\section{$\mathcal{SW}_1$ vs. $\mathcal{W}_1$}
\noindent \textbf{Example.} Consider $\mathbb{R}^3.$ We let $\mu$ be the uniform probability measure on the unit disc on the $xy$-plane and $\nu$ be the disc shifted up by $\varepsilon >0$ on the $z$-axis. 
\\
\\
Then, $\nu = T_{\#}\mu$ where $T(x) = x + k\varepsilon,$ where $k$ denotes the unit vector on the $z$-axis. Then the cost is 
\[
\int||T(x) - x||d\mu(x) = \int ||x+k\varepsilon -x||d\mu(x) = \int||k\varepsilon||d\mu(x) =\varepsilon.
\]
We don't know that $T$ is optimal, so $\mathcal{W}_1(\mu,\nu) \leq \varepsilon.$ To show the other direction, take $f(x,y,z) = z,$ which is 1-Lipschitz. For all $v \in\text{supp}(\mu),$ the $z$-coordinate of $v$ is 0. Thus $\int f d\mu = 0.$ Then for all $u \in \text{supp}(\nu),$ the $z$-coordinate of $u$ is $\varepsilon.$ So, 
\[
\int fd\nu - \int fd\mu = \varepsilon - 0 = \varepsilon,
\]
Then Kantorovich duality gives us $\mathcal{W}_1(\mu,\nu) \geq \varepsilon$ since $\mathcal{W}_1$ is the supremum over all such differences. So $\mathcal{W}_1(\mu,\nu) = \varepsilon.$
\\
\\
Now consider $\mathcal{SW}_1(\mu,\nu),$
\[
\mathcal{SW}_1(\mu,\nu) = \int_{\mathbb{S}^2} \mathcal{W}_1(P_{\theta\#}\mu, P_{\theta\#}\nu)d\sigma(\theta). 
\]
Notice that $P_{\theta\#}\nu =P_{\theta\#}T_\#\mu = (P_{\theta} \circ T)_{\#}\mu.$ Since $P_{\theta}$ is linear, $P_\theta(T(x)) = P_\theta(x) +\varepsilon P_{\theta}(k).$ Note $\varepsilon P_{\theta}(k) = \varepsilon \text{cos}(\phi)$ since $\theta = (\text{sin}(\phi) \text{cos}(\psi), \, \text{sin}(\phi)\text{sin}(\psi), \, \text{cos}(\phi))$ in spherical coordinates. Since $P_{\theta\#}\nu$ is  $P_{\theta\#}\mu$ shifted by $\varepsilon \text{cos}({\phi})$, the same argument as above gives
\[
\mathcal{W}_1(P_{\theta\#}\mu, P_{\theta\#}\nu) = \varepsilon |\text{cos}\phi|.
\]
Since $d\sigma$ is the uniform probability measure on $\mathbb{S}^2,$ we normalize so the integral evaluates to 1: $d\sigma = \frac{1}{4\pi}\text{sin}(\phi)d\phi d\psi$. Then,
\[
\mathcal{SW}_1(\mu,\nu) = \int_{\mathbb{S}^2} \varepsilon|\text{cos}(\phi)|d\sigma(\theta) = \frac{\varepsilon}{4\pi} \int_{0}^{2\pi}d\psi \int_{0}^{\pi}|\cos(\phi)|\sin(\phi)d\phi = \frac{\varepsilon}{2}.
\]
Therefore, 
\[
\frac{\mathcal{SW}_1}{\mathcal{W}_1} = \frac{1}{2}.
\]
In this example we notice that the ratio is independent of our choice of $\varepsilon$. The case is particularly simple, since $\nu$ is simply $\mu$ translated across the $z$-axis, and $\mathcal{SW}_1$ and $\mathcal{W}_1$ scale linearly with respect to $\varepsilon.$ Of course, other examples will not have such a simple structure. We proceed with a proposition and argument on translated measures for general dimension $d.$
\\
\\
\textbf{Proposition 1.}\textit{
Let $V$ be a $d$-dimensional real Hilbert space. Let $\mu \in \mathcal{P}_1(V)$ be a probability measure, and let $\nu =T_{\#}\mu$ where $T(x) = x + v, \ x, v \in V$ with $v \neq 0$. Then $\mathcal{W}_1(\mu,\nu) = ||v||.$}
\begin{proof}
Notice that by letting $y=T(x),$ 
\[
||x-T(x)|| = ||x - x - v|| = ||v||.
\]
So averaging over $\mu,$
\[
\int_{V} ||v||d\mu(x) = ||v|| .
\]
Thus, $\mathcal{W}_1(\mu,\nu) \leq ||v||.$ To show the other direction, let $u = \frac{v}{||v||}$ for $v \neq 0.$ Define $f(x) = \langle u, x \rangle.$ We claim that $f$ is 1-Lipschitz. To see, notice that for any $x,y,$ $|f(x) - f(y)| = |\langle u, x-y \rangle| \leq ||u||||x-y||,$ by Cauchy-Schwarz. Note $||u|| = 1,$ so $|f(x)-f(y)| \leq ||x-y||.$ Thus $f$ is 1-Lipschitz. Then by Kantorovich Duality,
\[
\mathcal{W}_1(\mu,\nu) = \sup_{\text{Lip(f)} \leq 1} \left( \int fd\nu - \int fd\mu\right).
\]
Using the fact that $\nu = T_{\#}\mu,$ 
\[
\int fd\nu = \int f(T(x))d\mu(x) = \int f(x+v)d\mu(x). 
\]
Then since $f(x) = \langle u, x\rangle,$
\[
f(x+v) = \langle u , x \rangle + \langle u, v \rangle = f(x) + ||v||.
\]
Since $\langle u, v \rangle  = \langle \frac{v}{||v||}, v  \rangle = \frac{||v||^2}{||v||} = ||v||.$ Therefore, 
\begin{align*}
\int fd\nu  = \int fd\mu + \int ||v||d\mu 
\\
= \int fd\mu + ||v|| 
\\\implies \int fd\nu - \int fd\mu = ||v||.
\end{align*}
since $\mu$ is a probability measure. Then we have that $||v|| \leq \mathcal{W}_1(\mu, \nu)$ since $\mathcal{W}_1$ is obtained by the supremum of all such 1-Lipschitz functions, so finally $\mathcal{W}_1(\mu, \nu) = ||v||.$
\end{proof}
Intuitively, if we obtain one measure $\nu$ by shifting another measure $\mu$ by a vector $v$, the total "work'' that Wasserstein distance captures is the length of the vector, $||v||.$ Indeed, in the previous example in $\mathbb{R}^3,$ $\mathcal{W}_1(\mu, \nu) = ||k\varepsilon|| = \varepsilon.$
\\
\\
\textbf{Proposition 2.} \textit{Let $V$ be a real $d$-dimensional Hilbert space. Let $\mu \in \mathcal{P}_1(V)$ be a probability measure, and let $\nu =T_{\#}\mu$ where $T(x) = x + v, \ x,v \in V$ with $v \neq 0.$ Then for $d \geq 2, $ $\mathcal{SW}_1(\mu, \nu) = \frac{||v||  2   \Gamma(\frac{d}{2})}{(d-1)\sqrt{\pi}\Gamma\big(\frac{d-1}{2}\big)}$ where $\Gamma$ denotes the Gamma function.}
\begin{proof}
Since $\nu =T_{\#}\mu,$
\[
P_{\theta\#}\nu = (P_{\theta} \circ T)_{\#}\mu.
\]
\\
Fix $x \in V.$ Then, $P_{\theta}(T(x)) = P_{\theta}(x+v) = P_{\theta}(x) + \langle v, \theta \rangle.$ So $P_{\theta \#} \nu$ is $P_{\theta \#}\mu $ shifted by $\langle v, \theta \rangle$ on $\mathbb{R}$. So by Proposition $1,$ 
\[
\mathcal{W}_1\bigl(P_{\theta\#}\mu, P_{\theta\#}\nu\bigr) = |\langle v, \theta \rangle| \ \forall\theta \in \mathbb{S}^{d-1}.
\]
Solving for $\mathcal{SW}_1(\mu,\nu)$ reduces to integrating the above quantity over $\theta \in \mathbb{S}^{d-1},$
\[
\mathcal{SW}_1(\mu,\nu) = \int_{\mathbb{S}^{d-1}} \mathcal{W}_1(P_{\theta\#}\mu, P_{\theta\#}\nu)d\sigma(\theta) = \int_{\mathbb{S}^{d-1}} |\langle v,\theta \rangle| d\sigma(\theta).
\]
Since the uniform measure $\sigma$ on $\mathbb{S}^{d-1}$ is rotation invariant for all $d,$ we choose coordinates so that $v$ is in the direction of the "north pole''. Let $e := v/|v|$ be the point where $\text{span}(v)$ intersects with $\mathbb{S}^{d-1}$ in the positive direction. Since any point $x \in \mathbb{S}^{d-1}$ can be written as [4, Eq. (20)]
\[
x = \pm [\cos \phi , u\sin \phi ]^T, \ \phi \in [0, \pi]
\]
where $u$ is a $d-1$ dimensional unit vector, we deconstruct $d\sigma(\theta)$ \cite[Eq. (21)]{Kent2022}:
\[
d\sigma(\theta) = \frac{\sin^{d-2}\phi d\phi du}{\int_0^{\pi}\sin^{d-2}\phi d\phi}. 
\] 
Notice that integrating both sides yields 
\[
\int_{\mathbb{S}^{d-1}}d\sigma(\theta) = \frac{\int_{0}^{\pi}\sin^{d-2}\phi d\phi \int_{\mathbb{S}^{d-2}}du}{\int_0^{\pi}\sin^{d-2}\phi d\phi} = \int_{\mathbb{S}^{d-2}}du = 1,
\]
since $\sigma$ is a probability measure. Plugging this back into the original integral, we get 
\begin{align*}
\int_{\mathbb{S}^{d-1}} |\langle v, \theta \rangle|d\sigma(\theta) 
\\
= ||v||\int_{\mathbb{S}^{d-1}}|\cos\phi|d\sigma(\theta) 
\\
= ||v|| \cdot \frac{\int_{0}^\pi |\cos\phi|\sin^{d-2} \phi d\phi}
{\int_0^\pi \sin^{d-2}\phi d\phi}.
\end{align*}
The numerator is an elementary integral solved by $u$-substitution. By symmetry,
\[
\int_0^{\pi} |\cos\phi|\sin^{d-2}\phi d\phi = 2\int_0^{\pi/2}|\cos\phi|\sin^{d-2}\phi d\phi.
\]
Then $\cos \phi \geq 0$ for $\phi \in [0,\pi/2]$ so we omit the absolute value. Let $u = \sin \phi, du = \cos \phi d\phi.$ Substituting,
\[
2 \int_0^{1} u ^{d-2}du = \frac{2}{d-1}.
\]
Using symmetry again, the integral in the denominator can be simplified as 
\[
\int_{0}^\pi \sin^{d-2} \phi d\phi = 2 \int_{0}^{\pi/2} \sin^{d-2} \phi d\phi 
\]
Then the integral is one of the well known Wallis integrals \cite[Eq. (5.6)]{Artin1964}, which has the form
\[
W_{d-2} =  \frac{\Gamma\big(\frac{d-1}{2}\big)\Gamma\big(\frac{1}{2}\big)}{2\Gamma\big(\frac{d}{2}\big)}.
\]
Now we see,
\[
2 \int_{0}^{\pi/2} \sin^{d-2} \phi d\phi = 2W_{d-2} = \frac{\Gamma\big(\frac{d-1}{2}\big)\Gamma\big(\frac{1}{2}\big)}{\Gamma\big(\frac{d}{2}\big)} = \frac{\Gamma\big(\frac{d-1}{2}\big)\sqrt{\pi}}{\Gamma\big(\frac{d}{2}\big)}.
\]
We finally obtain 
\[
\mathcal{SW}_1(\mu,\nu)=\frac{||v|| 2  \Gamma(\frac{d}{2})}{(d-1) \sqrt{\pi} \Gamma\big(\frac{d-1}{2}\big)}.
\]
\end{proof}
\noindent \textbf{Corollary 1.} \textit{$\frac{\mathcal{SW}_1}{\mathcal{W}_1}$ for translated measures.}
\[
\frac{\mathcal{SW}_1}{\mathcal{W}_1} = \frac{2   \Gamma(\frac{d}{2})}{(d-1)\sqrt{\pi}\Gamma\big(\frac{d-1}{2}\big)}.
\]
\\
\section{Lower bound on $c_d$}
We now provide an explicit lower bound for the constant $c_d$ for translated measures whose existence was established by Carlier, Figalli, Mérigot, and Wang \cite{CFMW2025}. Since their proof only establishes the existence of $c_d$, there is no explicit value currently given. Using Propositions 1 and 2 we construct such a lower bound for translated measures. Define:
\[
M_d = \frac{2\Gamma(\frac{d}{2})}{(d-1)\sqrt{\pi}\Gamma\big(\frac{d-1}{2}\big)}.
\]
\textbf{Theorem 1.} \textit{For any dimension $d \geq 2$ and all $R > 0$, let $\mu \in \mathcal{P}_1(\mathbb{R}^d)$ and $\nu = T_{\#}\mu$ where $T(x) = x + v$, $x, v \in \mathbb{R}^d$ with $v \neq 0$, and suppose $\mu, \nu$ are both supported on $B_R$. Then}
\[
M_d^{-\frac{1}{d}} \cdot \left(\frac{||v||}{R}\right)^{\frac{d-1}{d}} \leq c_d.
\]
\begin{proof}
    Let $d \geq 2$. Then by Carlier, Figalli, Mérigot, and Wang, there exists a constant $c_d$ such that for all $R > 0,$ 
\[
\mathcal{W}_1(\mu,\nu) \leq c_dR^{\frac{d-1}{d}}\mathcal{SW}_1(\mu,\nu)^{\frac{1}{d}}.
\]
\\
Then by Proposition $1$ and $2$ we substitute values for $\mathcal{W}_1$ and $\mathcal{SW}_1,$ 
\[
||v|| \leq c_dR^{\frac{d-1}{d}}\left(\frac{||v|| 2  \Gamma(\frac{d}{2})}{(d-1) \sqrt{\pi} \Gamma\big(\frac{d-1}{2}\big)}\right)^{1/d}. 
\]
Dividing both sides, the left hand side becomes 
\[
\frac{||v||}{||v||^{\frac{1}{d}}R^{\frac{d-1}{d}}} \cdot \left(\frac{2  \Gamma(\frac{d}{2})}{(d-1) \sqrt{\pi} \Gamma\big(\frac{d-1}{2}\big)}\right)^{-\frac{1}{d}} 
= \left(\frac{||v||}{R}\right)^{\frac{d-1}{d}} \cdot M_d^{-\frac{1}{d}}.
\]
Thus,
\[
M_d^{-\frac{1}{d}} \cdot \left(\frac{||v||}{R}\right)^{\frac{d-1}{d}} \leq c_d.
\]
\end{proof}
It is important to note that this bound comes from one specific family of measures, translations, and that the bound has not been proven to be tight. A natural direction for future work is computing the ratio $\frac{\mathcal{SW}_1}{\mathcal{W}_1}$ for other families of measures, such as rotations or scalings, and comparing various lower bounds on $c_d.$ 
\\
\\
\textbf{Acknowledgments.} LLMs (Claude Sonnet 4.6/5, ChatGPT Plus, Gemini) were used in the mathematical discussion and assisted with the editing of the manuscript. The author has independently verified all work and takes full responsibility for the contents of the paper.

\end{document}